\documentclass{amsart}%
\usepackage{amsfonts}
\usepackage[hiresbb]{graphicx}
\usepackage{amsmath}
\usepackage{amssymb}
\usepackage{version}%
\providecommand{\U}[1]{\protect\rule{.1in}{.1in}}
\theoremstyle{theorem}
\newtheorem{defn}{Definition}[section]
\newtheorem{thm}[defn]{Theorem}
\newtheorem{exmp}[defn]{Example}

\begin{document}
\title[Generating Functions of Periodic Recurrences]{On the generating functions of a periodic infinite order linear recurrence}
\author[A. Bravo \& H. M. Oliveira]{Ant\'{o}nio Bravo$^{a}$ and Henrique M. Oliveira$^{b}$}
\address{$^{a}$Centro de An\'{a}lise Funcional e Aplica\c{c}\~{o}es\\
$^{b}$Centro de An\'{a}lise Matem\'{a}tica Geometria e Sistemas Din\^{a}micos\\
$^{a,b}$Mathematics Department, Technical Institute of Lisbon, University of Lisbon\\
Av. Rovisco Pais, 1049-001 Lisbon, Portugal\\
Email Adresses: $^{a}$abravo@math.ist.utl.pt, $^{b}$holiv@math.ist.utl.pt
(corresponding author)}
\date{November 4th, 2013 -- AMS: 15A15, 39A06, 92D25,39A06}
\maketitle

\begin{abstract}
Let $\mathcal{G}$ be the space of generating functions of a periodic infinite
order linear recurrence. In this paper we provide an explicit procedure for
computing a basis of $\mathcal{G}$.

\end{abstract}

\section{Introduction and Statement of the Main Result}

Let $\mathbb{N}=\left\{  0,1,2,\ldots\right\}  $ be the set of the natural
numbers and $p$ a positive integer. In this paper we study a particular class
of vectorial homogeneous recurrences of the type%

\begin{equation}
\mathrm{x}_{n+1}=\sum\limits_{i=0}^{+\infty}\mathbf{A}_{i}\mathrm{x}%
_{n-i}\text{, for all }n\in\mathbb{N}\text{,} \label{vector-rec}%
\end{equation}
where $\mathrm{x}_{n}\in\mathbb{C}^{p}$ and $\left\{  \mathbf{A}_{n}\right\}
_{n\in\mathbb{N}}$ is an infinite sequence of $p\times p$ matrices with
complex entries.

As in \cite{rash}, we call this type of infinite order linear recurrences
\emph{generalized Fibonacci recurrences}, or $Fib_{p}^{\infty}$ recurrences
for short. A $Fib_{p}^{\infty}$ recurrence is said to be $s$-periodic (with
$s\in\mathbb{Z}^{+}$) if $\mathbf{A}_{n+s}=\mathbf{A}_{n}$ for all
$n\in\mathbb{N}.$

The main result of this work concerns the study of the solutions of a
$s$-periodic $Fib_{p}^{\infty}$ recurrence.

Before stating the main result of this work, we introduce some notation.

Throughout the paper, $\mathbb{C}\left[  \left[  z\right]  \right]  $ denotes
the commutative ring of formal power series with complex coefficients.
Matrices with entries in $\mathbb{C}$ and $\mathbb{C}\left[  \left[  z\right]
\right]  $ will be denoted respectively as elements of $\mathbb{C}^{m\times
n}$ and $\mathbb{C}\left[  \left[  z\right]  \right]  ^{m\times n}$. As usual,
the entry $\left(  i,j\right)  $ of a matrix $\mathbf{A}$ will be denoted by
$\mathbf{A}\left(  i,j\right)  \mathbf{.}$

The infinite-dimensional vector spaces over $\mathbb{C}$%
\[
U=\bigoplus\limits_{n\in\mathbb{N}}\mathbb{C}^{p}\text{ and }V=%
{\textstyle\prod\limits_{n\in\mathbb{Z}}}
\mathbb{C}^{p}\text{,}%
\]
will play an important role in this discussion. We write $\mathbf{u}$ and
$\mathbf{v}$ to denote the vectors of $U$ and $V$, with components
$\mathrm{u}_{n}\in\mathbb{C}^{p}$ and $\mathrm{v}_{n}\in\mathbb{C}^{p}$, i.e.,%
\[
\mathbf{u}=\left(  \mathrm{u}_{n}\right)  _{n\in\mathbb{N}}=\left(
\mathrm{u}_{0},\mathrm{u}_{1},\ldots\right)  \text{, with }\mathrm{u}_{n}%
\in\mathbb{C}^{p}\text{, }\mathrm{u}_{n}=0\text{ for almost all }n\text{,}%
\]
and%
\[
\mathbf{v}=\left(  \mathrm{v}_{n}\right)  _{n\in\mathbb{N}}=\left(
\ldots,\mathrm{v}_{-1},\mathrm{v}_{0},\mathrm{v}_{1},\ldots\right)  \text{,
with }\mathrm{v}_{n}\in\mathbb{C}^{p}\text{.}%
\]
We write $\mathrm{e}_{1},\ldots,\mathrm{e}_{p}$ to denote the vectors of the
standard basis of $\mathbb{C}^{p}$. The symbols $\mathbf{e}_{\beta}$, with
$\beta\in\mathbb{Z}^{+}$, denote the vectors of the standard basis of $U$:%
\[%
\begin{array}
[c]{c}%
\mathbf{e}_{1}=\left(  \mathrm{e}_{1},\mathrm{0},\mathrm{0},\ldots\right)
\text{, }\mathbf{e}_{2}=\left(  \mathrm{e}_{2},\mathrm{0},\mathrm{0}%
,\ldots\right)  \text{, }\ldots\text{, }\mathbf{e}_{p}=\left(  \mathrm{e}%
_{p},\mathrm{0},\mathrm{0},\ldots\right)  \text{,}\\
\mathbf{e}_{p+1}=\left(  \mathrm{0},\mathrm{e}_{1},\mathrm{0},\ldots\right)
\text{, }\mathbf{e}_{p+2}=\left(  \mathrm{0},\mathrm{e}_{2},\mathrm{0}%
,\ldots\right)  \text{, }\ldots\text{, }\mathbf{e}_{2p}=\left(  \mathrm{0}%
,\mathrm{e}_{p},\mathrm{0},\ldots\right)  \text{,}\\
\ldots
\end{array}
\]
where $\mathrm{0}$ denotes the zero vector of $\mathbb{C}^{p}$.

A vector $\mathbf{v}=\left(  \mathrm{v}_{n}\right)  _{n\in\mathbb{N}}$ $\in V$
is said to be a solution of the $Fib_{p}^{\infty}$ recurrence
\eqref{vector-rec}, if the set $\left\{  n\leq0:\mathrm{v}_{n}\neq
\mathrm{0}\right\}  $ is finite and
\[
\mathrm{v}_{n+1}=\sum\limits_{i=0}^{+\infty}\mathbf{A}_{i}\mathrm{v}%
_{n-i}\text{, for all }n\in\mathbb{N}\text{.}%
\]
The subspace of $V$ whose vectors are the solutions of the $Fib_{p}^{\infty}$
recurrence is denoted by $S$.

Naturally, there exists an isomorphism $\Theta:U\mapsto S$, where
\[
\Theta\left(  \mathbf{u}\right)  =\left(  v_{n}\right)  _{n\in\mathbb{N}%
}\text{,}%
\]
denotes the unique vector of $S$ satisfying%
\[
v_{-n}=u_{n}\text{ for all }n\in\mathbb{N}\text{.}%
\]
The vector $\Theta\left(  \mathbf{u}\right)  \in S$ is called the solution of
the $Fib_{p}^{\infty}$ recurrence for the initial condition $\mathbf{u\in}U$.
The vector space $U$ is called space of initial conditions.

In order to analyze the asymptotic behaviour of a solution $\Theta\left(
\mathbf{u}\right)  =\left(  v_{n}\right)  _{n\in\mathbb{Z}}\in S$, we define
the generating function $G\left(  \mathbf{u}\right)  $ as the formal power
series with coefficients in $\mathbb{C}^{p}$%
\[
G\left(  \mathbf{u}\right)  =\sum\limits_{n\geq0}v_{n}z^{n}\text{.}%
\]
Alternatively, $G\left(  \mathbf{u}\right)  $ can be defined as an element of
the $\mathbb{C}$-vector space $\mathbb{C}\left[  \left[  z\right]  \right]
^{p}$ such that%
\[
G\left(  \mathbf{u}\right)  =\left(  G_{1}\left(  \mathbf{u}\right)
,\ldots,G_{p}\left(  \mathbf{u}\right)  \right)  \text{,}%
\]
with%
\[
G_{\alpha}\left(  \mathbf{u}\right)  =\sum\limits_{n\geq0}v_{n}^{\left(
\alpha\right)  }z^{n}\in\mathbb{C}\left[  \left[  z\right]  \right]  \text{,}%
\]
where $v_{n}^{\left(  \alpha\right)  }$ denotes the $\alpha$-component of
$v_{n}$ with respect to the standard base of $\mathbb{C}^{p}$.

Naturally, the map
\[%
\begin{array}
[c]{cccc}%
G & :U & \mapsto & \mathbb{C}\left[  \left[  z\right]  \right]  ^{p}\\
& \mathbf{u} & \rightarrow & G\left(  \mathbf{u}\right)  \text{,}%
\end{array}
\]
is linear. We call $\mathcal{G}=G\left(  U\right)  $ the space of generating
functions of the $Fib_{p}^{\infty}$ recurrence.

Notice that, from the linearity of the map $G$, we have%

\begin{equation}
G\left(  \mathbf{u}\right)  =\sum\limits_{\beta\geq1}c_{\beta}G\left(
\mathbf{e}_{\beta}\right)  \text{,} \label{f2}%
\end{equation}
where $\left(  c_{\beta}\right)  _{\beta\in\mathbb{Z}^{+}}$ denote the
coordinates of $\mathbf{u}$ with respect to the standard basis, $\left(
\mathbf{e}_{\beta}\right)  _{\beta\in\mathbb{Z}^{+}}$, of $U$. Therefore, the
set%
\[
\left\{  G\left(  \mathbf{e}_{\beta}\right)  :\beta\in\mathbb{Z}^{+}\right\}
\text{,}%
\]
spans $\mathcal{G}$.

Now, we can state the main result of this work concerning the important class
of periodic $Fib_{p}^{\infty}$ recurrences. As we will see, this result
combined with the main theorem of \cite{jalves1} provides an explicit
procedure for computing a basis of $\mathcal{G}$.

\begin{thm}
\label{t1}If a $Fib_{p}^{\infty}$ recurrence is $s$-periodic, then all
generating functions $G\left(  \mathbf{u}\right)  \in\mathcal{G}$ are rational
function of $z$. Moreover, the generating functions $G\left(  \mathbf{e}%
_{1}\right)  $ $,\ldots,$ $G\left(  \mathbf{e}_{\left(  s+1\right)  p}\right)
$ span the space $\mathcal{G}$. Hence, $\dim\mathcal{G}\leq\left(  s+1\right)
p$.
\end{thm}

\section{Proof of Theorem \ref{t1}}

The notions of kneading matrix and kneading determinant of a $Fib_{p}^{\infty
}$ recurrence, introduced in \cite{jalves1} (see also \cite{jalves2}), are the
main ingredients of the proof of Theorem \ref{t1}. In order to improve the
readability of the proof of Theorem \ref{t1}, we present a brief description
of the main ideas of \cite{jalves1}.

Let $\alpha=1,\ldots,p$ and $\beta$ be a positive integer. Analogously to the
kneading theory, stated by Milnor and Thurston in \cite{knead}, we introduce
the $\left(  \alpha,\beta\right)  $-kneading increment of a $Fib_{p}^{\infty}$
recurrence as the following formal power series%
\[
K\left(  \alpha,\beta\right)  =\left\{
\begin{array}
[c]{ll}%
\sum_{n\geq0}\mathbf{A}_{n+q-1}\left(  \alpha,p\right)  z^{n}\in
\mathbb{C}\left[  \left[  z\right]  \right]  \text{,} & \text{if }p\text{
divides }\beta\text{,}\\
\sum_{n\geq0}\mathbf{A}_{n+q}\left(  \alpha,r\right)  z^{n}\in\mathbb{C}%
\left[  \left[  z\right]  \right]  \text{,} & \text{otherwise},
\end{array}
\right.
\]
where $q$ and $r$ denote, respectively, the quotient and the remainder of the
integer division of $\beta$ by $p$.

Notice that, in the particular case $\alpha,\beta=1,\ldots,p$ the kneading
increment is given by%
\[
K\left(  \alpha,\beta\right)  =\underset{n\geq0}{\sum}\mathbf{A}_{n}\left(
\alpha,\beta\right)  z^{n}\text{,}%
\]
and the $p\times p$-matrix of formal power series%
\[
\mathbf{K}=\left(
\begin{array}
[c]{ccc}%
K\left(  1,1\right)  & \cdots & K\left(  1,p\right) \\
\vdots & \ddots & \vdots\\
K\left(  p,1\right)  & \cdots & K\left(  p,p\right)
\end{array}
\right)  \in\mathbb{C}\left[  \left[  z\right]  \right]  ^{p\times p}\text{,}%
\]
receives the name of kneading matrix of the $Fib_{p}^{\infty}$ recurrence.

Finally, we define kneading determinant of the $Fib_{p}^{\infty}$ recurrence
by%
\[
\Delta=\det\left(  \mathbf{I}-z\mathbf{K}\right)  \in\mathbb{C}\left[  \left[
z\right]  \right]  \text{,}%
\]
where $\mathbf{I}$ denotes the $p\times p$-identity matrix.

Similarly, for each $\alpha=1,\ldots,p$ and $\beta\in\mathbb{Z}^{+}$ we define
the extended kneading matrix $\mathbf{K}_{\alpha}\left(  \beta\right)  $ and
the extended kneading determinant $\Delta_{\alpha}\left(  \beta\right)  $ of a
$Fib_{p}^{\infty}$ recurrence by setting%
\begin{equation}
\mathbf{K}_{\alpha}\left(  \beta\right)  =\left(
\begin{array}
[c]{cccc}%
K\left(  1,1\right)  & \cdots & K\left(  1,p\right)  & K\left(  1,\beta\right)
\\
\vdots & \ddots & \vdots & \vdots\\
K\left(  p,1\right)  & \cdots & K\left(  p,p\right)  & K\left(  p,\beta\right)
\\
\delta\left(  \alpha,1\right)  & \cdots & \delta\left(  \alpha,p\right)  &
\delta\left(  \alpha,\beta\right)
\end{array}
\right)  \in\mathbb{C}\left[  \left[  z\right]  \right]  ^{\left(  p+1\right)
\times\left(  p+1\right)  }\text{,} \label{f4}%
\end{equation}
where $\delta\left(  i,j\right)  $ is the usual Kronecker delta function and%
\[
\Delta_{\alpha}\left(  \beta\right)  =\det\left(  \mathbf{I}-z\mathbf{K}%
_{\alpha}\left(  \beta\right)  \right)  \in\mathbb{C}\left[  \left[  z\right]
\right]  \text{,}%
\]
where $\mathbf{I}$ denotes the $\left(  p+1\right)  \times\left(  p+1\right)
$ identity matrix.

Next, we recall the main result of \cite{jalves1}, which establishes a main
relationship between the generating functions $G_{\alpha}\left(
\mathbf{e}_{\beta}\right)  $ and the kneading determinants $\Delta$ and
$\Delta_{\alpha}\left(  \beta\right)  $ of a $Fib_{p}^{\infty}$ recurrence.

\begin{thm}
\label{T2}For every $\alpha=1,..,p$ and every vector $\mathbf{e}_{\beta}$ of
the standard basis of $U$, the generating function $G_{\alpha}\left(
\mathbf{e}_{\beta}\right)  $ of a $Fib_{p}^{\infty}$ recurrence satisfies the
following equation in $\mathbb{C}\left[  \left[  z\right]  \right]  $%
\[
zG_{\alpha}\left(  \mathbf{e}_{\beta}\right)  =1-\Delta^{-1}\Delta_{\alpha
}\left(  \beta\right)  \text{.}%
\]

\end{thm}

We have all the ingredients needed to prove Theorem \ref{t1}.

First, observe that for any $\alpha=1,\ldots,p$ and $\beta\in\mathbb{Z}^{+}$
the kneading increment $K\left(  \alpha,\beta\right)  $ of a $s$-periodic
$Fib_{p}^{\infty}$ recurrence is a rational function of $z$. Indeed, as
$\mathbf{A}_{n+s}=\mathbf{A}_{n}$ for all $n\in\mathbb{N}$, one gets%
\begin{align}
K\left(  \alpha,\beta\right)   &  =\left\{
\begin{array}
[c]{ll}%
\sum_{n\geq0}\mathbf{A}_{n+q-1}\left(  \alpha,p\right)  z^{n}\text{,}%
\vspace{0.15cm} & \text{if }p\text{ divides }\beta\text{,}\\
\sum_{n\geq0}\mathbf{A}_{n+q}\left(  \alpha,r\right)  z^{n}\text{,} &
\text{otherwise,}%
\end{array}
\right. \label{f3}\\
&  =\left\{
\begin{array}
[c]{ll}%
\sum_{m\geq0}\left(  \sum_{k=0}^{s-1}\mathbf{A}_{k+q-1}\left(  \alpha
,p\right)  z^{k}\right)  z^{sm}\text{,}\vspace{0.15cm} & \text{if }p\text{
divides }\beta\text{,}\\
\sum_{m\geq0}\left(  \sum_{k=0}^{s-1}\mathbf{A}_{k+q}\left(  \alpha,r\right)
z^{k}\right)  z^{sm}\text{,} & \text{otherwise,}%
\end{array}
\right. \nonumber\\
&  =\left\{
\begin{array}
[c]{ll}%
\left(  \sum_{k=0}^{s-1}\mathbf{A}_{k+q-1}\left(  \alpha,p\right)
z^{k}\right)  \sum_{m\geq0}z^{sm}\text{,}\vspace{0.15cm} & \text{if }p\text{
divides }\beta\text{,}\\
\left(  \sum_{k=0}^{s-1}\mathbf{A}_{k+q}\left(  \alpha,r\right)  z^{k}\right)
\sum_{m\geq0}z^{sm}\text{,} & \text{otherwise,}%
\end{array}
\right. \nonumber\\
&  =\left\{
\begin{array}
[c]{ll}%
\frac{\sum_{k=0}^{s-1}\mathbf{A}_{k+q-1}\left(  \alpha,p\right)  z^{k}%
}{1-z^{s}}\text{,}\vspace{0.15cm} & \text{if }p\text{ divides }\beta\text{,}\\
\frac{\sum_{k=0}^{s-1}\mathbf{A}_{k+q}\left(  \alpha,r\right)  z^{k}}{1-z^{s}%
}\text{,} & \text{otherwise,}%
\end{array}
\right. \nonumber
\end{align}
where $q$ and $r$ denote, respectively, the quotient and the remainder of the
integer division of $\beta$ by $p$. So, the kneading determinant $\Delta$, as
well any extended kneading determinant $\Delta_{\alpha}\left(  \beta\right)  $
of a $s$-periodic $Fib_{p}^{\infty}$ recurrence are rational functions of $z$.
But, by Theorem \ref{T2}, this means that any generation function $G_{\alpha
}\left(  \mathbf{e}_{\beta}\right)  $ is rational too. Combining this with
\eqref{f2} we finally conclude that any generating function $G\left(
\mathbf{u}\right)  $ of a $s$-periodic $Fib_{p}^{\infty}$ recurrence is a
rational function of $z$. This is precisely what is stated in the first part
of Theorem \ref{t1}.

Next, we prove the second part of Theorem \ref{t1}.

Evidently, if $q\in\mathbb{N}$ and $r\in\mathbb{N}$ denote, respectively, the
quotient and the remainder of the division of $\beta\in\mathbb{Z}^{+}$ by $p$,
then the quotient and the remainder of the division of $\beta+sp$ by $p$ are,
respectively, $q+s$ and $r$. So, since the $Fib_{p}^{\infty}$ recurrence is
$s$-periodic, one gets by \eqref{f3}%
\begin{align*}
K\left(  \alpha,\beta+sp\right)   &  =\left\{
\begin{array}
[c]{ll}%
\frac{\sum_{k=0}^{s-1}\mathbf{A}_{k+q+s-1}\left(  \alpha,p\right)  z^{k}%
}{1-z^{s}}\text{,}\vspace{0.15cm} & \text{if }p\text{ divides }\beta
+sp\text{,}\\
\frac{\sum_{k=0}^{s-1}\mathbf{A}_{k+q+s}\left(  \alpha,r\right)  z^{k}%
}{1-z^{s}}\text{,} & \text{otherwise,}%
\end{array}
\right. \\
&  =\left\{
\begin{array}
[c]{ll}%
\frac{\sum_{k=0}^{s-1}\mathbf{A}_{k+q-1}\left(  \alpha,p\right)  z^{k}%
}{1-z^{s}}\text{,}\vspace{0.15cm} & \text{if }p\text{ divides }\beta
+sp\text{,}\\
\frac{\sum_{k=0}^{s-1}\mathbf{A}_{k+q}\left(  \alpha,r\right)  z^{k}}{1-z^{s}%
}\text{,} & \text{otherwise,}%
\end{array}
\right. \\
&  =\left\{
\begin{array}
[c]{ll}%
\frac{\sum_{k=0}^{s-1}\mathbf{A}_{k+q-1}\left(  \alpha,p\right)  z^{k}%
}{1-z^{s}}\text{,}\vspace{0.15cm} & \text{if }p\text{ divides }\beta\text{,}\\
\frac{\sum_{k=0}^{s-1}\mathbf{A}_{k+q}\left(  \alpha,r\right)  z^{k}}{1-z^{s}%
}\text{,} & \text{otherwise,}%
\end{array}
\right. \\
&  =K\left(  \alpha,\beta\right)  \text{,}%
\end{align*}
for all $\alpha=1,\ldots,p$ and $\beta\in\mathbb{Z}^{+}$. Combining this with
\eqref{f4}, one obtains%
\[
\mathbf{K}_{\alpha}\left(  \beta\right)  =\mathbf{K}_{\alpha}\left(
\beta+sp\right)  \text{, for all }\alpha=1,\ldots,p\text{ and }\beta>p\text{.}%
\]
Therefore,%
\[
\Delta_{\alpha}\left(  \beta\right)  =\det\left(  \mathbf{I}-z\mathbf{K}%
_{\alpha}\left(  \beta\right)  \right)  =\det\left(  \mathbf{I}-z\mathbf{K}%
_{\alpha}\left(  \beta+sp\right)  \right)  =\Delta_{\alpha}\left(
\beta+sp\right)  \text{,}%
\]
and, by Theorem \ref{T2}, we can write%
\[
G_{\alpha}\left(  \mathbf{e}_{\beta}\right)  =G_{\alpha}\left(  \mathbf{e}%
_{\beta+sp}\right)  \text{, for all }\alpha=1,\ldots,p\text{ and }%
\beta>p\text{.}%
\]
This proves that $G\left(  \mathbf{e}_{\beta+sp}\right)  =G\left(
\mathbf{e}_{\beta}\right)  $ for all $\beta>p$. Consequently,
\[
\left\{  G\left(  \mathbf{e}_{\beta}\right)  :\beta\in\mathbb{Z}^{+}\right\}
=\left\{  G\left(  \mathbf{e}_{\beta}\right)  :\beta=1,\ldots,p\left(
1+s\right)  \right\}  \text{,}%
\]
which proves that $G\left(  \mathbf{e}_{1}\right)  ,G\left(  \mathbf{e}%
_{2}\right)  ,\ldots,G\left(  \mathbf{e}_{p\left(  1+s\right)  }\right)  $
span $\mathcal{G}$.

Hence, $\dim\mathcal{G}\leq p\left(  1+s\right)  $. This is precisely what is
stated in the second part of Theorem \ref{t1}.

We finish the paper with two simple examples to illustrate the role played by
Theorems \ref{t1} and \ref{T2} in the explicit computation of a basis of
$\mathcal{G}$.

\begin{exmp}
Consider the $1$-periodic $Fib_{2}^{\infty}$ recurrence%
\[
\mathrm{x}_{n+1}=\sum\limits_{i=0}^{+\infty}\mathbf{A}_{i}\mathrm{x}%
_{n-i}\text{, for all }n\in\mathbb{N}\text{,}%
\]
with%
\[
\mathbf{A}_{n}\mathbf{=}\left(
\begin{array}
[c]{cc}%
1 & 2\\
3 & 4
\end{array}
\right)  \text{, for all }n\in\mathbb{N}.
\]
By \eqref{f3}, the kneading increments are%
\[
K\left(  1,\beta\right)  =\left\{
\begin{array}
[c]{ll}%
\frac{2}{1-z}\text{,}\vspace{0.15cm} & \text{if }\beta\text{ is even,}\\
\frac{1}{1-z}\text{,} & \text{if }\beta\text{ is odd,}%
\end{array}
\right.  \text{ and }K\left(  2,\beta\right)  =\left\{
\begin{array}
[c]{ll}%
\frac{4}{1-z}\text{,}\vspace{0.15cm} & \text{if }\beta\text{ is even,}\\
\frac{3}{1-z}\text{,} & \text{if }\beta\text{ is odd,}%
\end{array}
\right.
\]
for all $\beta\in\mathbb{Z}^{+}$. Therefore,%
\begin{equation}
\mathbf{K}=\left(
\begin{array}
[c]{cc}%
\frac{1}{1-z} & \frac{2}{1-z}\\
\frac{3}{1-z} & \frac{4}{1-z}%
\end{array}
\right)  \text{, }\Delta=\det\left(  \mathbf{I}-z\mathbf{K}\right)
=\frac{4z^{2}-7z+1}{\left(  1-z\right)  ^{2}}\text{.} \label{f5}%
\end{equation}
On the other hand, as%
\[
\mathbf{K}_{\alpha}\left(  \beta\right)  =\left\{
\begin{array}
[c]{cc}%
\left(
\begin{array}
[c]{ccc}%
\frac{1}{1-z} & \frac{2}{1-z} & \frac{2}{1-z}\\
\frac{3}{1-z} & \frac{4}{1-z} & \frac{4}{1-z}\\
\delta\left(  \alpha,1\right)  & \delta\left(  \alpha,2\right)  &
\delta\left(  \alpha,\beta\right)
\end{array}
\right)  \text{,}\vspace{0.15cm} & \text{if }\beta\text{ is even,}\\
\left(
\begin{array}
[c]{ccc}%
\frac{1}{1-z} & \frac{2}{1-z} & \frac{1}{1-z}\\
\frac{3}{1-z} & \frac{4}{1-z} & \frac{3}{1-z}\\
\delta\left(  \alpha,1\right)  & \delta\left(  \alpha,2\right)  &
\delta\left(  \alpha,\beta\right)
\end{array}
\right)  \text{,} & \text{if }\beta\text{ is odd,}%
\end{array}
\right.
\]
One has $\Delta_{\alpha}\left(  \beta\right)  =\det\left(  \mathbf{I}%
-z\mathbf{K}_{\alpha}\left(  \beta\right)  \right)  =$%
\begin{equation}
\left\{
\begin{array}
[c]{l}%
\begin{array}
[c]{c}%
\frac{\delta\left(  \alpha,\beta\right)  \left(  7z^{2}-4z^{3}-z\right)
-2\delta\left(  \alpha,2\right)  z^{2}+2z^{3}+2z^{2}-7z+1}{\left(  1-z\right)
^{2}}%
\end{array}
\text{, if }\beta\text{ is even,}\\%
\begin{array}
[c]{c}%
\frac{\delta\left(  \alpha,\beta\right)  \left(  -4z^{3}+7z^{2}-z\right)
+2\delta\left(  \alpha,2\right)  z^{2}\left(  2z-1\right)  -z^{3}+3z^{2}%
-7z+1}{\left(  1-z\right)  ^{2}}%
\end{array}
\text{, if }\beta\text{ is odd.}%
\end{array}
\right.  \label{f6}%
\end{equation}
From Theorem \ref{t1}, the generating functions $G\left(  \mathbf{e}%
_{1}\right)  $, $G\left(  \mathbf{e}_{2}\right)  $, $G\left(  \mathbf{e}%
_{3}\right)  $ and $G\left(  \mathbf{e}_{4}\right)  $ span $\mathcal{G}$. By
\eqref{f5}, \eqref{f6} and Theorem \ref{T2}, the computation of the generating
functions is straightforward
\[
G\left(  \mathbf{e}_{1}\right)  =\left(  \frac{5z^{2}-6z+1}{4z^{2}-7z+1}%
,\frac{-3z^{2}+3z}{4z^{2}-7z+1}\right)  \text{,}%
\]%
\[
G\left(  \mathbf{e}_{2}\right)  =\left(  \frac{-2z^{2}+2z}{4z^{2}-7z+1}%
,\frac{2z^{2}-3z+1}{4z^{2}-7z+1}\right)  \text{,}%
\]%
\[
G\left(  \mathbf{e}_{3}\right)  =\left(  \frac{z^{2}+z}{4z^{2}-7z+1}%
,\frac{-3z^{2}+3z}{4z^{2}-7z+1}\right)  \text{,}%
\]%
\[
G\left(  \mathbf{e}_{4}\right)  =\left(  \frac{-2z^{2}+2z}{4z^{2}-7z+1}%
,\frac{-2z^{2}+4z}{4z^{2}-7z+1}\right)  \text{.}%
\]
Thus, $G\left(  \mathbf{e}_{1}\right)  $, $G\left(  \mathbf{e}_{2}\right)  $,
$G\left(  \mathbf{e}_{3}\right)  $, $G\left(  \mathbf{e}_{4}\right)  $ are
linearly independent. Therefore, $\dim$ $\mathcal{G}=4$.
\end{exmp}

\begin{exmp}
Now, consider the $1$-periodic $Fib_{2}^{\infty}$ recurrence%
\[
\mathrm{x}_{n+1}=\sum\limits_{i=0}^{+\infty}\mathbf{A}_{i}\mathrm{x}%
_{n-i}\text{, for all }n\in\mathbb{N}\text{,}%
\]
with%
\[
\mathbf{A}_{n}=\left(
\begin{array}
[c]{cc}%
1 & 1\\
1 & 1
\end{array}
\right)  \text{, for all }n\in\mathbb{N}\text{.}%
\]
As in the previous example, the generating functions $G\left(  \mathbf{e}%
_{1}\right)  $, $G\left(  \mathbf{e}_{2}\right)  $, $G\left(  \mathbf{e}%
_{3}\right)  $, $G\left(  \mathbf{e}_{4}\right)  $ span $\mathcal{G}$. Now, by
Theorem \ref{T2}, one gets
\[
G\left(  \mathbf{e}_{1}\right)  =\left(  \frac{2z-1}{3z-1},\frac{-z}%
{3z-1}\right)  \text{, }G\left(  \mathbf{e}_{2}\right)  =\left(  \frac
{-z}{3z-1},\frac{2z-1}{3z-1}\right)  \text{,}%
\]
and%
\[
G\left(  \mathbf{e}_{3}\right)  =G\left(  \mathbf{e}_{4}\right)  =\left(
\frac{-z}{3z-1},\frac{-z}{3z-1}\right)  \text{.}%
\]
Thus, $G\left(  \mathbf{e}_{1}\right)  $, $G\left(  \mathbf{e}_{2}\right)  $,
$G\left(  \mathbf{e}_{3}\right)  $ are linearly independent and $\dim
\mathcal{G}=3$.
\end{exmp}

\section*{Acknowledgements}

The authors were supported by CAMGSD-LARSys and CEAF through FCT Program POCTI-FEDER.

\bibliographystyle{plain}
\bibliography{biblioH}

\label{lastpage}
\end{document}